\documentclass[11pt]{amsart}
\usepackage{amsmath,amssymb,amsfonts,url,mathptmx}
\usepackage{color}
\usepackage{appendix}
\numberwithin{equation}{section}

\newtheorem{thm}{Theorem}[section]

\usepackage{amsthm,a0size,bm,indentfirst,wasysym}

\usepackage{amsmath}

\allowdisplaybreaks[4]

\begin{document}
	\title[Simple blow-up]{Simple blow-up solutions of singular Liouville equations}
	\keywords{Liouville equation, quantized singular source, non-simple blow-up, construction of solutions, blow-up solutions. Spherical Harnack inequality}

\author{Lina Wu}\footnote{Lina Wu is partially supported by National Natural Science Foundation of China (12201030), China Postdoctoral Science Foundation (2022M720394) and Talent Fund of Beijing Jiaotong University (2022RC028). }

	\address{  Lina Wu\\
		School of Mathematics and Statistics \\
		Beijing Jiaotong University \\
		Beijing, 100044, China }
	\email{lnwu@bjtu.edu.cn}

	\date{\today}
	
	\begin{abstract}
		In a recent series of important works \cite{wei-zhang-1,wei-zhang-2,wei-zhang-3}, Wei-Zhang proved several vanishing theorems for non-simple blow-up solutions of singular Liouville equations. It is well known that a non-simple blow-up situation happens when the spherical Harnack inequality is violated near a quantized singular source. In this article, we further strengthen the conclusions of Wei-Zhang by proving that if the spherical Harnack inequality does hold, there exist blow-up solutions with non-vanishing coefficient functions.
	\end{abstract}
	
	
	\maketitle
	
\section{Introduction}

It is well known that the following Liouville equation has a rich background in geometry and Physics. 
\begin{equation}\label{liou-eq}
\Delta u+h(x)e^{u(x)}=\sum_{t=1}^L4\pi \gamma_t \delta_{p_t}
\quad \mbox{in}\quad \Omega\subset \mathbb R^2, 
\end{equation}
where $\Omega$ is a subset of $\mathbb R^2$, $p_1,..p_L$ are $L$ points in $\Omega$ and $4\pi\gamma_t\delta_{p_t}$ ($t=1,...,L$) are Dirac masses placed at $p_t$. Since applications require integrability of $e^{u}$ we assume $\gamma_t>-1$ for each $t$. 

Equation \eqref{liou-eq} is one of the most extensively studied elliptic partial differential equations in recent years. In conformal geometry, \eqref{liou-eq} is related to the well-known Nirenberg problem when all $\gamma_t=0$. The recent progress on this project can be seen in Kazdan-Warner \cite{kazdan-warner}, Chang-Gursky-Yang \cite{chang-gursky-yang}, Chang-Yang \cite{chang-yang}, Cheng-Lin \cite{cheng-lin}, and the references therein. If some $\gamma_t\neq 0$, \eqref{liou-eq} arises from the existence of conformal metric with conic singularities, seen in Fang-Lai \cite{fang-lai}, Troyanov \cite{troy}, Wei-Zhang \cite{wei-zhang-pacific}. Also, it serves as
a model equation in the Chern-Simons-Higgs theory and in the Liouville system, the interested readers may browse Chanillo-Kiessling \cite{chani-kiess}, Spruck-Yang \cite{spruck-yang}, Tarantello \cite{taran-C}, Yang \cite{yang-y}, and the references therein.

It is well known that if there is no singularity in (\ref{liou-eq}), $h\equiv 1$ and $\int_{\mathbb R^2}e^u{\rm d}x<\infty$, a global solution belongs to a family described by three parameters (see \cite{chen-li-duke}). Then Y. Y. Li \cite{li-cmp} proved the first uniform approximation theorem, which confirms that around a regular blow-up point, the profile of a blow-up sequence is close to that of a sequence of global solutions. Later  Chen-Lin \cite{chen-lin-sharp}, Zhang \cite{zhang-cmp}, Gluck \cite{gluck}, Bartolucci, et,al \cite{unique-4}
improved Li's estimate by obtaining better pointwise estimates and some gradient estimates. It turns out that the blow-up point has to be a critical point of a function determined by the coefficient function. This plays a crucial role in applications. In the non-quantized case, the classification theorem was proved by Prajapat-Tarantello, the uniform estimate is obtained by Bartolucci-Chen-Lin-Tarantello \cite{bart-chen-lin-tarant}, Bartolucci-Tarantello \cite{bart-taran}, Zhang \cite{zhang-ccm}. The most difficult case is when the singular source is quantized. In this case, the first breakthrough was obtained by Kuo-Lin in \cite{kuo-lin}, then independently by Bartolucci-Tarantello in \cite{bart-taran}. In this case, if the spherical Harnack inequality is violated near a quantized singular source, the profile of bubbling solutions appears to have multiple local maximums. Here a sequence of bubbling solutions satisfying \emph{spherical Harnack inequality} means the oscillation of solutions on each fixed radius around the singular point is uniformly bounded. In the work of Kuo-Lin, they use \emph{non-simple blow-up} to describe this phenomenon. In a recent series of works of Wei-Zhang \cite{wei-zhang-1,wei-zhang-2,wei-zhang-3}, they proved the first vanishing theorems for the non-simple blow-up case. Their two main results can be stated as follows: 

Let $\{u_k\}_{k=1}^{\infty}$ be a sequence of blow-up solutions of 
 \begin{equation}\label{main-2}
        \Delta u_k+|x|^{2\alpha}h_k(x)e^{u_k(x)}=0, \quad \mbox{in}\quad B_1 
    \end{equation}
where $h_k$ is a sequence of smooth, positive functions in $B_1$:
\begin{equation}\label{h-assump}
	\frac 1{c_1}\le h_k(x)\le c_1, \quad \|\nabla^{\beta}h_k(x)\|_{B_1}\le c_1, \quad x\in B_1, \quad |\beta |=1,2,3. 
	\end{equation}
for some $c_1>0$. Let $0$ be the only blow-up point of $u_k$ in $B_1$,
	and suppose $u_k$ has a bounded oscillation on $\partial B_1$:
 \begin{equation}\label{BOF}
 |u_k(x)-u_k(y)|\le C,\quad \forall x,y\in \partial B_1,
 \end{equation}
 and a uniform bound on its integration:
 \begin{equation}\label{unif-en}
 \int_{B_1}|x|^{2\alpha}h_k(x)e^{u_k(x)}{\rm d}x<C
 \end{equation}
 for some $C>0$ independent of $k$. In their first vanishing theorem Wei-Zhang proved that 

 \emph{Theorem A: (Wei-Zhang). Let $u_k$ be a sequence of non-simple blow-up solutions around the origin. Suppose $0$ is the only blow-up point in $B_1$ and
 $u_k$ satisfies (\ref{main-2}),(\ref{BOF}) and (\ref{unif-en}). Then along a sub-sequence
 $$\lim_{k\to \infty}\nabla (\log h_k+\psi_k)(0)=0$$ 
 where $\psi_k$ is the harmonic function that eliminates the finite oscillation of $u_k$ on $\partial B_1$:
 \begin{equation}\label{psi-eq}
 \Delta \psi_k=0,\quad \mbox{in}\quad B_1, \quad \psi_k(x)=u_k(x)-\frac{1}{2\pi}\int_{\partial B_1}u_k,\quad x\in \partial B_1. 
 \end{equation}
 }

 In their recent work, Wei-Zhang further proved the following Laplacian vanishing theorem:

 \emph{Theorem B: (Wei-Zhang). Let $u_k$ be the same as in Theorem A. Then along a subsequence, 
 $$\lim_{k\to \infty}\Delta (\log h_k)(0)=0.$$
 }

 It is important to point out that in both Theorem A and Theorem B, the blow-up sequence has to be \emph{non-simple}, this assumption implies that $\alpha\in \mathbb N$. Both Theorem A and Theorem B are powerful tools in application, since the equation (\ref{main-2}) represents a number of situations in more general equations/systems. For example, in the author's recent joint work with Wei and Zhang \cite{wei-wu-zhang}, we proved that under certain conditions on the coefficient function and Gauss curvature, all blow-up points to Toda systems are simple. 

The purpose of this article is twofold. First if $\alpha=0$ in (\ref{main-2}) and $0$ is the only blow-up point, it is well known that (see \cite{chen-lin-sharp, gluck,zhang-cmp}) along a sub-sequence
$\lim_{k\to \infty} \nabla (\log h_k+\psi_k)(0)=0$. Over the years it has long been suspected that this property does not hold if $\alpha$ is not an integer. This is indeed verified in our first main theorem: 
\begin{thm}\label{main-thm} 
For any given $\alpha>-1$ and $\alpha\not\in \mathbb N\cup \{0\}$, there exist a sequence $h_k$ satisfying (\ref{h-assump}) and
\begin{equation}\label{h-assum-2}
	|\nabla \log h_k(0)+\nabla \psi_k(0)|\ge c_1,\quad |\Delta \log h_k(0)|\ge c_1
	\end{equation}
	for some $c_1>0$,
 Corresponding to $h_k$ there is a sequence of blow-up solutions $u_k$ of (\ref{main-2}) such that the origin is its only blow-up point, (\ref{BOF}) (\ref{unif-en}) holds for $u_k$, which also satisfies the spherical Harnack inequality around the origin. 
\end{thm}

The second goal is to prove that 
when $\alpha\in \mathbb N\cup \{0\}$ we can construct a sequence of \emph{simple } blow-up solutions that
does not satisfy the Laplacian vanishing theorem. 

\begin{thm}\label{main-thm-2} 
Let $\alpha\in \mathbb N\cup \{0\}$, there exist a sequence of blow-up solutions $\{u_k\}_{k=1}^{\infty}$ of (\ref{main-2}) having $0$ as its only blow-up point in $B_1$. Moreover $\{u_k\}$ satisfies (\ref{BOF}) (\ref{unif-en}) and the coefficient $h_k$ satisfies (\ref{h-assump}) and 
$$|\Delta (\log h_k)(0)|\ge c,  \,\, \mbox{for a constant }\,\, c>0 \,\, \mbox{independent of $k$}.$$
\end{thm}

 Theorem \ref{main-thm} settles the conjecture that around a non-quantized singular source, the vanishing theorems do not hold. Theorem \ref{main-thm-2} proves that it is essential to have a \emph{ non-simple} blow-up sequence in Theorem B. If this assumption is violated, the corresponding Laplacian vanishing property also fails. However, this article did not provide a similar example for the first-order vanishing theorem in Theorem A. 

 The paper is organized as follows: In Section \ref{non-quan}, we establish Theorem \ref{main-thm}. Our proof is based on the thorough comprehension of the corresponding linearized operator of a model equation. It is also essential that we analyze the Fourier series of some correction terms and prove its convergence. In Section \ref{quan}, we establish Theorem \ref{main-thm-2}, and the key point of the proof is to use a radial coefficient function and reduce all the iterations into radial cases. This method made us avoid kernel functions in the linearized equation corresponding to the quantized case.

\section{Non-quantized situation}\label{non-quan}
	
	In this section, we consider the non-quantized case. In other words,  we set $\alpha>-1$ and $\alpha\notin\mathbb{N}\cup\{0\}$. It is known that the spherical Harnack holds around the origin when $\alpha$ is not an integer (See \cite{kuo-lin}).
	
	Denote $\lambda_k=u_k(0)$ and $\epsilon_k=e^{-\frac{\lambda_k}{2(1+\alpha)}}$. Let $v_k$ be the scaling of $u_k$: 
	\begin{equation*}
		v_k(y)=u_k(\epsilon_ky)+2(1+\alpha)\log \epsilon_k,\quad y\in\Omega_k:=B(0,\epsilon_k^{-1}).
	\end{equation*}
	Clearly, we need to construct $v_k$ to satisfy
	\begin{equation}\label{euq-v-k}
		\left\{\begin{array}{lcl}
			\Delta v_k(y)+|y|^{2\alpha}h_k(\epsilon_ky)e^{v_k(y)}=0,&& {\rm in} \ \, \Omega_k,	\\
			v_k(0)=0,   \\
			|v_k(y_1)-v_k(y_2)|\le C,&& {\rm for\ \ any} \ \, y_1,y_2\in\partial\Omega_k,\\
			v_k(y)\to-2\log(1+|y|^{2+2\alpha}),&& {
				\rm in}\ \, C_{loc}^{\beta}(\mathbb{R}^2)
		\end{array}
		\right.
	\end{equation}
	where $\beta\in (0,1)$.
	It suffices to construct $\{v_k\}$ satisfying \eqref{euq-v-k}. Since we can choose $h_k$ we require $h_k(0)= 8(1+\alpha)^2$ for convenience. Let
	\begin{equation*}
		U_k(y)=-2\log(1+|y|^{2+2\alpha})
	\end{equation*}
	be a standard bubble that satisfies
	\begin{equation}\label{equ-U-k}
		\Delta U_k(y)+8(1+\alpha)^2|y|^{2\alpha}e^{U_k(y)}=0 \quad {\rm in} \ \, \mathbb{R}^2.
	\end{equation}
	Here we note that a uniform estimate of Bartolucci-Chen-Lin-Tarantello \cite{bart-chen-lin-tarant} assures that any blow-up solution $v_k$ of (\ref{euq-v-k}) satisfies 
	\begin{equation*}
		|v_k(y)-U_k(y)|\leq C,\quad y\in \Omega_k.
	\end{equation*}

	We will construct our solutions based on the expansion of $v_k$ established in \cite{zhang-ccm}. Firstly, let us recall some notations and results in \cite{zhang-ccm}. Denote
	\begin{equation*}
		\begin{split}
			g_k(r)&=-\frac{1}{4\alpha(1+\alpha)}\frac{r}{1+r^{2+2\alpha}},\quad r=|y|, \\
			c_1^k(y)&=g_k(r)\epsilon_k\sum_{j=1}^{2}\partial_jh_k(0)\theta_j,\quad \theta_j=\frac{y_j}{r}\ \ (j=1,2).
		\end{split}
	\end{equation*}
	Then $c_1^k$ satisfies 
	\begin{equation}\label{equ-phi-k}
		\Delta c_1^k+8(1+\alpha)^2|y|^{2\alpha}e^{U_k(y)}c_1^k=-\sum_{j=1}^{2}\epsilon_k\partial_jh_k(0)y_j|y|^{2\alpha}e^{U_k(y)}\quad {\rm in} \ \, \Omega_k.
	\end{equation}
	\cite{zhang-ccm} tells us that $c_1^k$ is the second term in the expansion of $v_k$ if $\alpha>0$ is a non-integer. For the case $-1<\alpha<0$, Bartolucci-Yang-Zhang \cite{byz} have established the same result.
 Here we point out that the radial part of $c_1^k$ decays like $\epsilon_k r^{-1-2\alpha}$ at infinity. In particular, for $r=\epsilon_k^{-1}$, the angular part of the function is comparable to $\epsilon_k^{2+2\alpha}e^{i\theta}$, which means this term contributes no oscillation on the boundary. So as long as $|\nabla \log h_k(0)|\ge 2c>0$, we have $|\nabla \log h_k(0)+\nabla \psi_k(0)|\ge c$. 

For the convenience of the readers, we comment that the construction of $c_1^k$ is essentially solving
\begin{equation}\label{equ-l}
    \frac{d^2}{dr^2}g+\frac{1}{r}\frac{d}{dr}g+\Big(8(1+\alpha)^2r^{2\alpha}e^{U_k}-\frac{l^2}{r^2}\Big)g=-r^{1+2\alpha}e^{U_k},\quad r>0,
\end{equation}
with $l=1$. From the proof of Lemma 2.1 in \cite{zhang-ccm}, we know two fundamental solutions $F_1$ and $F_2$ of the homogeneous equation of \eqref{equ-l} can be written explicitly as follows:
    \begin{equation}\label{two-fun}
	    \begin{split}
	        F_{1}(r)&=\frac{(\frac{l}{1+\alpha}+1)r^l+(\frac{l}{1+\alpha}-1)r^{l+2(1+\alpha)}}{1+r^{2(1+\alpha)}}, \\
	        F_{2}(r)&=\frac{(\frac{l}{1+\alpha}+1)r^{-l+2(1+\alpha)}+(\frac{l}{1+\alpha}-1)r^{-l}}{1+r^{2(1+\alpha)}}.
	    \end{split}
	\end{equation}
Therefore, we can verify that $g$ can be explicitly written with two fundamental solutions above by the standard ODE methods. 

 For the motivation of adding more terms in the correction, we use the decay of $c_1^k$ to
 obtain
 \begin{align}\label{tem-2}
     &\Delta (U_k+c_1^k)+8(1+\alpha)^2r^{2\alpha}e^{U_k+c_1^k}\\
     =&\Delta U_k+\Delta c_1^k+8(1+\alpha)^2r^{2\alpha}e^{U_k}\big(1+c_1^k+\frac{(c_1^k)^2}{2}\big)
    +O(\epsilon_k^3)(1+r)^{-7-8\alpha}\nonumber \\
    =&8(1+\alpha)^2r^{2\alpha}e^{U_k}\frac{(c_1^k)^2}{2}-\epsilon_k\sum_j\partial_jh_k(0)\theta_jr^{1+2\alpha}e^{U_k}+O(\epsilon_k^3)(1+r)^{-7-8\alpha}.\nonumber
 \end{align}
 
  At this moment we write the expansion of $h_k(\epsilon_ky)$:
	\begin{equation}\label{h-expan}
        \begin{split}
            h_k(\epsilon_ky)=&8(1+\alpha)^2+\epsilon_k\nabla h_k(0)\cdot y+\frac{\epsilon_k^2}2\partial_{11}h_k(0)(y_1^2-\frac{|y|^2}{2})
	+\frac{\epsilon_k^2}2\partial_{22}h_k(0)\cdot \\
 &(y_2^2-\frac{|y|^2}{2}) 
	+\epsilon_k^2\partial_{12}h_k(0)
	y_1y_2+\frac{\epsilon_k^2}4\Delta h_k(0)|y|^2+O(\epsilon_k^3)|y|^3  \\
 =&8(1+\alpha)^2+\epsilon_k\nabla h_k(0)\cdot y+\epsilon_k^2r^2\Theta_2+\frac 14\epsilon_k^2 r^2\Delta h_k(0)+O(\epsilon_k^3)r^3.
        \end{split}
	\end{equation}
 where
 \begin{equation*}
    \begin{split}
        \Theta_2:=&\frac 12\partial_{11}h_k(0)(\theta_1^2-\frac 12)+\partial_{12}h_k(0)\theta_1\theta_2+\frac 12\partial_{22}h_k(0)(\theta_2^2-\frac 12)\\
     =&\frac 14(\partial_{11}h_k(0)-\partial_{22}h_k(0))\cos 2\theta +\frac 12 \partial_{12}h_k(0)\sin 2\theta. 
    \end{split}
 \end{equation*}
 Based on (\ref{equ-phi-k}), (\ref{tem-2}) and (\ref{h-expan}) we have
 \begin{align}\label{euq-Uk-c1k}
 \begin{split}
     &\Delta (U_k+c_1^k)+h_k(\epsilon_ky)|y|^{2\alpha}e^{U_k+c_1^k}\\
 =&r^{2\alpha}e^{U_k}\Big(8(1+\alpha)^2\frac{(c_1^k)^2}2+\Theta_2\epsilon_k^2r^2+\frac 14 \epsilon_k^2r^2\Delta h_k(0)+\epsilon_k\nabla h_k(0)\cdot y c_1^k \\
 &+O(\epsilon_k^3 (1+r)^3)\Big).
 \end{split}
 \end{align}
 
 Now we compute $(c_1^k)^2$:
 \begin{align*}
     (c_1^k)^2=&\epsilon_k^2g_k^2(\partial_1 h_k(0) \cos \theta+\partial_2 h_k(0) \sin\theta)^2\\
     =&\epsilon_k^2g_k^2\Big(\frac{|\nabla h_k(0) |^2}2+\frac 12
     ((\partial_1 h_k(0))^2-(\partial_2h_k(0))^2)\cos 2\theta  \\
     &+\partial_1 h_k(0)\partial_2 h_k(0) \sin 2\theta\Big)
     \end{align*}
     Also the remaining term of the order $O(\epsilon_k^2)$ is 
\begin{align*}
    &\epsilon_k\nabla h_k(0)\cdot y c_1^k\\
 =&\epsilon_k^2 g_k r (\partial_1h_k(0)\cos \theta+\partial_2 h_k(0)\sin \theta)^2\\
 =&\epsilon_k^2 g_k r \Big(\frac{|\nabla h_k(0)|^2}2+\frac 12(\partial_1h_k(0)^2-\partial_2 h_k(0)^2)\cos 2\theta+\partial_1h_k(0)\partial_2 h_k(0)\sin 2\theta\Big).
\end{align*}

 To get rid of the terms with $e^{2i\theta}$ of the order $O(\epsilon_k^2)$ in \eqref{euq-Uk-c1k} we let $c_2^k$ be the solution of
 \begin{align*}
 &\Delta c_2^k+8(1+\alpha)^2|y|^{2\alpha}e^{U_k}c_2^k=-r^{2\alpha}e^{U_k}\Big(\epsilon_k^2r^2\Theta_2+\mathcal{A}(\frac{(c_1^k)^2}{2}+\epsilon_k\nabla h_k(0)\cdot yc_1^k)\Big)\\
 =&-\epsilon_k^2|y|^{2\alpha}e^{U_k}\Big (r^2\Theta_2
 +\big(4(1+\alpha)^2g_k^2+g_kr\big)\cdot 
 ((\partial_1h_k(0))^2-(\partial_2 h_k(0))^2)\frac{\cos 2\theta}2\\
 &+\partial_1h_k(0) \partial_2 h_k(0)\sin 2\theta)\Big ). 
 \end{align*}
Note that $\mathcal{A}(\cdot)$ means the non-radial part of the term in the parenthesis. 

Since each term in $c_2^k$ is a product of a radial function and a spherical harmonic function, we set $w_1$ to be a solution of 
$$\frac{d^2}{dr^2}w_1+\frac{1}{r}\frac{d}{dr}w_1+\Big(8(1+\alpha)^2r^{2\alpha}e^{U_k}-\frac{4}{r^2}\Big)w_1=r^{2+2\alpha}e^{U_k}$$
with the control of $|w_1(r)|\le C$ for all $r$. Similarly, we set $w_2$ to be a solution of 
$$\frac{d^2}{dr^2}w_2+\frac{1}{r}\frac{d}{dr}w_2+\Big(8(1+\alpha)^2r^{2\alpha}e^{U_k}-\frac{4}{r^2}\Big)w_2=r^{2\alpha}e^{U_k}\big(4(1+\alpha)^2g_k^2+g_kr\big)$$
with $|w_2(r)|\le C$ for all $r$. Two fundamental solutions of the corresponding homogeneous equation can be seen in \eqref{two-fun} with $l=2$. Furthermore, we observe that the non-homogenous terms have good decay rates at infinity. Therefore, the construction of $w_1$ and $w_2$ is standard. At this point, it is easy to verify that $c_2^k$ can be constructed as 
\begin{align*}
c_2^k(y)
=&\epsilon_k^2\Big (w_1(r)\Theta_2+w_2(r) \big ((\partial_1h_k(0))^2-(\partial_2 h_k(0))^2)\frac{\cos 2\theta}2 \\
&+\partial_1h_k(0) \partial_2 h_k(0)\sin 2\theta\big )\Big).
\end{align*}

 Finally we use $c_0^k$ to handle the radial term of the order $O(\epsilon_k^2)$:
We let $c_0^k$ solve
 \begin{align*}
&\Delta c_0^k+8(1+\alpha)^2|y|^{2\alpha}e^{U_k}c_0^k\\
=&-\epsilon_k^2|y|^{2\alpha}e^{U_k}\Big (r^2\frac{\Delta h_k(0)}4+\frac 12|\nabla h_k(0)|^2\big(4(1+\alpha)^2g_k(r)^2+g_k(r)r\big)\Big ).
 \end{align*}
Since both $U_k$ and the right-hand side of the above are radial, we can construct $c_0^k$ as a radial function $c_0^k(r)$ that satisfies
$$\left\{\begin{array}{ll}
\frac{d^2}{dr^2}c_0^k(r)+\frac{d}{dr}c_0^k(r)+8(1+\alpha)^2r^{2\alpha}e^{U_k}c_0^k(r)\\
\quad =-\epsilon_k^2r^{2\alpha}e^{U_k}\Big (\frac{r^2}4\Delta h_k(0)+\frac 12|\nabla h_k(0)|^2\big(4(1+\alpha)^2g_k(r)^2+g_k(r)r\big)\Big ).\\
c_0^k(0)=\frac{d}{dr}c_0^k(0)=0.
\end{array}
\right.
$$
We only need to define $c_0^k$ for $0<r<\epsilon_k^{-1}$. It is easy to use the standard ODE method to obtain
\begin{equation}\label{est-c0}
|c_0^k(r)|\le C\epsilon_k^2(1+r)^{-2\alpha}\log (2+r),
\quad 0<r<\epsilon_k^{-1}.
\end{equation}
Set $c_k=c_0^k+c_1^k+c_2^k$, we verify by direct computation that 
\begin{equation}\label{base-e}
\Delta (U_k+c_k)+|y|^{2\alpha}h_k(\epsilon_ky)e^{U_k+c_k}
= E_k
\end{equation}
\begin{equation}\label{itera-2}
|E_k(y)|\le c_1\epsilon_k^3(1+|y|)^{-1-2\alpha},\quad y\in \Omega_k.
\end{equation}
So in order to find a solution with a non-vanishing coefficient, we need to find $d_k$ to satisfy
\begin{equation}\label{want-1}\Delta (U_k+c_k+d_k)+|y|^{2\alpha}h_k(\epsilon_ky)e^{U_k+c_k+d_k}=0,\quad \mbox{in}
\quad \Omega_k. 
\end{equation}

The difference between (\ref{base-e}) and (\ref{want-1}) gives
\begin{equation}\label{want-2}
\Delta d_k+8(1+\alpha)^2|y|^{2\alpha}e^{U_k}d_k=-E_k-f(d_k). 
\end{equation}
where 
\begin{equation}\label{itera-3}
f(d_k)=-|y|^{2\alpha}h_k(\epsilon_ky)e^{U_k+c_k}(e^{d_k}-1-d_k)
+|y|^{2\alpha}e^{U_k}(h_k(\epsilon_ky)e^{c_k}-h_k(0))d_k.
\end{equation}
 is of higher order. 
Based on (\ref{want-2}) we design an iteration scheme: Let $d_k^{(0)}\equiv 0$ and
$d_k^{(1)}$ satisfy
$$\Delta d_k^{(1)}+8(1+\alpha)^2|y|^{2\alpha}e^{U_k}d_k^{(1)}=-E_k-f(d_k^{(0)}). $$
In general we shall construct $d_k^{(m+1)}$ that satisfies
$$\Delta d_k^{(m+1)}+8(1+\alpha)^2|y|^{2\alpha}e^{U_k}d_k^{(m+1)}=-E_k-f(d_k^{(m)}) $$
and
\begin{equation}\label{der-d}
d_k^{(m+1)}(0)=|\nabla d_k^{(m+1)}(0)|=0.
\end{equation}
 
Here we claim that there exists $c_0>0$ independent of $m$ and $k$ such that 
\begin{equation}\label{unif-b}
|d_k^{(m)}(y)|\le c_0\epsilon_k^3(1+|y|)^{1-2\alpha}\log (2+|y|). 
\end{equation}
The constant $c_0$ will be determined based on $c_1$ later. 
To prove this uniform bound, we assume that (\ref{unif-b}) holds for $d_k^{(m)}$, and we shall show that it also holds for $d_k^{(m+1)}$. 

	The projection to $1$ is the following equation: Let $f_0$ be the 
 projection of $d_k^{(m+1)}$ onto $1$, then $f_0$ solves 
	\begin{equation*}
		\left\{\begin{array}{ll}
			f_{0}''(r)+\frac 1r f_{0}'(r)+8(1+\alpha)^2r^{2\alpha}e^{U_k}f_{0}(r)+E_0^k=0, \quad 0<r\leq\epsilon_k^{-1}\\
			f_{0}(0)=\frac{\mathrm{d}}{\mathrm{d}r}f_{0}(0)=0,
		\end{array}
		\right.
	\end{equation*}
 where $E_0^k$ is the corresponding projection of $E_k+f(d_k^{(m)})$ onto 1, and satisfies a similar bound of $E_k$:
 \begin{equation}\label{fur-err}
 |E_0^k(r)|\le 2c_1\epsilon_k^3(1+|y|)^{-1-2\alpha}.
 \end{equation}
 The reason that $E_0^k$ has a worse coefficient $2c_1$ is that the $d_k^{(m)}$ terms are absorbed. 
 
    We denote the two fundamental solutions of the homogeneous equation of $f_{0}$ as $u_1$ and $u_2$, where
	$$
	u_1(r)=\frac{1-r^{2+2\alpha}}{1+r^{2+2\alpha}},
	$$
	and $u_2(r)$ is comparable to $\log r$ near $0$ and infinity. Based on standard ODE theory, 
	\begin{equation*}
	    f_{0}(r)=-u_{1}(r)\int_0^{r} t E_1^k(t)u_{2}(t)\mathrm{d}t-u_{2}(r)\int_0^rt E_1^k(t)u_{1}(t)\mathrm{d}t.
	\end{equation*}
	Integrating the identity above, we know that
	\begin{equation*}
	    f_{0}(r)=O(\epsilon_k^3(1+r)^{1-2\alpha}r|\log(r)|),\ \ \mathrm{at}\ \ 0,\quad f_{0}(r)=O(\epsilon_k^3(1+r)^{1-2\alpha}\log(r)),\ \ \mathrm{at}\ \ \infty.
	\end{equation*}
	In other words, we have the following estimate for $f_{0}$
    $$ 
    |f_{0}(r)|\le  c_0\epsilon_k^3(1+r)^{1-2\alpha}\log(2+r),
    $$
    where $c_0$ is a constant independent of $l$ and only depends on $c_1$.

    Next, we consider the projections  on high frequencies. For $l\in\mathbb{N}^+$, let $f_{l}$ satisfies
	\begin{equation*}
	    \left\{\begin{array}{ll}
			f_{l}''(r)+\frac{1}{r}f_{l}'(r)+\Big(8(1+\alpha)^2r^{2\alpha}e^{U_k}-\frac{l^2}{r^2}\Big)f_{l}(r)+E_{2,k}^l=0, \quad 0<r\leq\epsilon_k^{-1}\\
			f_{l}(0)=0.
		\end{array}
		\right.
	\end{equation*}
	Here $E_{2,k}^l$ ($l\ge 1$) is the radial part of the projection of some error term on $\cos({l\theta}$): 
	\begin{equation*}
	    E_{2,k}^l(r)=\frac{1}{2\pi}\int_0^{2\pi}E_{k}\cos{(l\theta)}\mathrm{d}\theta
	\end{equation*}
	The estimate of $E_{2,k}^l$ is
	\begin{equation}\label{e2k-l}
	    |E_{2,k}^l(r)|\le 2 c_1\epsilon_k^3(1+r)^{-1-2\alpha}
	\end{equation}
    In order to find $f_{l}$ we use two fundamental solutions $F_{1}$ and $F_{2}$ of the homogeneous equation, whose explicit expressions can be seen in \eqref{two-fun}. 
    As one can see that $F_{1}$ is comparable to $r^l$ at the origin and at infinity, and $F_{2}$ is comparable to $r^{-l}$ at the origin and infinity. At this point,	we can construct $f_{l}$ as follows
	\begin{equation*}
	    f_{l}(r)=-F_{1}(r)\int_r^{\infty}\frac{t}{2l}E_{2,k}^l(t)F_{2}(t)\mathrm{d}t-F_{2}(r)\int_0^r\frac{t}{2l}E_{2,k}^l(t)F_{1}(t)\mathrm{d}t.
	\end{equation*}
    Integrate the identity above, we know that
    $$ 
    |f_{l}(r)|\le \frac{c_2}{l^2}\epsilon_k^3(1+r)^{1-2\alpha},
    $$
    where $c_2$ is a constant independent of $l$. It is easy to see that $f_{l}(0)=0$. Furthermore the summation of projections on all $\cos{(l\theta})$ $(l\ge 1)$ is convergent. That is
    \begin{equation*}
        \big|\sum_{l\ge 1}f_{l}(r)\big|\le c_2\epsilon_k^3(1+r)^{1-2\alpha}\sum_{l\ge1}\frac{1}{l^2}\le c_2\epsilon_k^3(1+r)^{1-2\alpha}
    \end{equation*}
    
    In the same way we can construct the projection  on $\sin{(l\theta)}$ for all $l\ge 1$, called $\tilde{f}_{l}$, and the summation of $\tilde{f}_{l}$ is convergent as well. $d_k^{(m)}$ is well-defined and satisfies the estimate \eqref{unif-b}.

Thus by Brower fixed point theorem, we obtain the existence of $d_k$. The construction is complete in this case. 

The Laplacian term is also obviously true, which can be seen in the construction. The construction of a non-quantized case is complete.

\section{Quantized situation}\label{quan}

Let $N$ be a positive natural number, our goal is to construct a sequence of blow-up solutions $u_k$ such that 
$$\Delta u_k+|x|^{2N}h_k(x)e^{u_k(x)}=0, \quad \mbox{in}\quad B_1 $$
such that the spherical Harnack holds around the origin, the only blow-up point in $B_1$ and $\Delta \log h_k(0)$ do not tend to zero. Here $\psi_k$ is the set of harmonic functions that eliminate the oscillation of $u_k$ on $\partial B_1$.

The main result of this section is to prove the following theorem.
\begin{thm}\label{lap-non}
For any $N\in \mathbb N$, there exists $h_k(x)$ satisfying \eqref{h-assump} and a sequence of blow-up solutions $u_k$ of \eqref{main-2}\eqref{BOF}\eqref{unif-en} such that $u_k$ is simple and $|\Delta (\log h_k)(0)|\ge c$ for some $c>0$ independent of $k$. 
\end{thm}

\noindent{\bf Proof of Theorem \ref{lap-non}:} We set 
$$h_k(x)=8(N+1)^2+|x|^2,\quad x\in B_1.$$
Obviously 
$$\nabla \log h_k(0)=0,\quad \Delta (\log h_k)(0)=\frac{\Delta h_k(0)}{h_k(0)}=\frac{1}{4(1+N)^2}.$$

Let $v_k$ be the scaling of $u_k$ according to the maximum of $u_k$: Let 
$$\epsilon_k=e^{-\frac{u_k(0)}{2(1+N)}}$$ and
$$v_k(y)=u_k(\epsilon_ky)+2(1+N)\log \epsilon_k. $$
The equation for $v_k$ is 
$$\Delta v_k+(1+\epsilon_k^2|y|^2)|y|^{2N}e^{v_k}=0.$$

Our goal is to construct $v_k$ satisfying the equation above based on the global solution $U_k$. The classification theorem of Prajapat-Tarantello \cite{prajapat} gives the standard bubble of $\Delta U+8(N+1)^{2}|y|^{2N}e^{U}=0$:
\begin{equation*}
    U(y)=\log \frac{\lambda}{\big(1+\lambda|y^{N+1}-\xi|^2\big)^2}
\end{equation*}
where parameters $\lambda>0$ and $\xi\in\mathbb{C}$.

Setting $\lambda=1$ and $\xi=0$ in $U$, we use the radial $U_k(y)$:
$$U_k(y)=\log \frac{1}{\big(1+|y|^{2N+2}\big)^2}. $$
Here we note that $\partial_{\lambda}|_{\lambda=0}
U$, $\partial_{\xi}|_{\xi=0} U$ and $\partial_{\bar \xi}|_{\xi=0}U$ for a basis for the linearized space.
\begin{align*}
    \partial_{\lambda}\big|_{\lambda=1}U&=\frac{1-r^{2N+2}}{1+r^{2N+2}}, \\
    \partial_{\xi}\big|_{\xi=0}U&=\frac{2r^{N+1}}{1+r^{2N+2}}e^{-i(N+1)\theta}, \\
    \partial_{\bar \xi}\big|_{\xi=0}U&=\frac{2r^{N+1}}{1+r^{2N+2}}e^{i(N+1)\theta}.
\end{align*}

Because of this, we see that corresponding to $N$ we have 
$$\frac{2r^{N+1}}{1+r^{2N+2}}\sin ((N+1)\theta),\quad
\frac{2r^{N+1}}{1+r^{2N+2}}\cos ((N+1)\theta) $$  in the kernel, this is the reason we only obtain the non-vanishing estimate for $\Delta (\log h_k)(0)$. It would be interesting to construct a simple blowup sequence with non-vanishing first-order coefficients. 

Based on the fact $h_k(\epsilon_ky)=8(N+1)^2+\epsilon_k^2|y|^2$ and the equation of $U_k$, we have
\begin{equation*}
    \Delta U_k+h_k(\epsilon_ky)|y|^{2N}e^{U_k}=\epsilon_k^2|y|^{2N+2}e^{U_k}.
\end{equation*}
In order to deal with the right-hand side of the equation above, we let $c_k$ solve  
\begin{equation*}
    \Delta c_k+8(N+1)^2|y|^{2N}e^{U_k}c_k=-\epsilon_k^2|y|^{2N+2}e^{U_k}.
\end{equation*}
Similar with $c_0^k$ in the non-quantized case, we can construct $c_k$ as a radial function $c_k(r)$ satisfying 
$$\left\{\begin{array}{ll}
\frac{d^2}{dr^2}c_k(r)+\frac{d}{dr}c_k(r)+8(N+1)^2r^{2N}e^{U_k}c_k(r)=-\epsilon_k^2r^{2N+2}e^{U_k},  \quad  0<r<\epsilon_k^{-1}\\
c_k(0)=\frac{d}{dr}c_k(0)=0.
\end{array}
\right.
$$
After the standard ODE method, we obtain the estimate as in \eqref{est-c0}:
\begin{equation}\label{est-c0-q}
|c_k(r)|\le C\epsilon_k^2(1+r)^{-2N}\log (2+r),
\quad 0<r<\epsilon_k^{-1}.
\end{equation}
Note that $e^{U_k+c_k}=e^{U_k}(1+c_k+O(\epsilon_k^4))$. By direct computation, we obtain 
\begin{equation}\label{base-e-q}
\Delta (U_k+c_k)+|y|^{2N}h_k(\epsilon_ky)e^{U_k+c_k}
= E_k.
\end{equation}
Here $E_k$ is radial and satisfies
\begin{equation}\label{itera-2-q}
|E_k(y)|\le c\epsilon_k^3(1+|y|)^{-1-2N},\quad y\in \Omega_k,
\end{equation}
where $c$ is a positive constant independent of $k$.

Then we set $v_k=U_k+c_k+b_k$. Removing the equation for $U_k$ and $c_k$ we write the equation of $b_k$ as
\begin{equation}\label{equ-bk}
\Delta b_k+|y|^{2N}h_k(\epsilon_ky)e^{U_k+c_k+b_k}-|y|^{2N}h_k(\epsilon_ky)e^{U_k+c_k}=-E_k,\quad |y|\le \tau \epsilon_k^{-1}. 
\end{equation}
The equation can be further written as
\begin{equation}\label{equ-bk-2}
\Delta b_k+8(N+1)^2|y|^{2N}e^{U_k}b_k=-E_k-\hat{f}(b_k),\quad |y|\le \tau \epsilon_k^{-1}
\end{equation}
where 
$$\hat{f}(b_k)=-|y|^{2N}h_k(\epsilon_k y)e^{U_k+c_k}(e^{b_k}-1-b_k)+|y|^{2N}e^{U_k}(h_k(\epsilon_ky)e^{c_k}-h_k(0))b_k.$$
Similar to the non-quantized case, we construct $b_k$ by iteration. Let $b_k^{(m)}\equiv 0$ and $b_k^{(1)}$ satisfy
\begin{equation*}
    \Delta b_k^{(1)}+8(N+1)^2r^{2N}e^{U_k}b_k^{(1)}=-E_k-\hat f(b_k^{(0)}).
\end{equation*}
In general we construct $b_k^{(m+1)}$ satisfying
\begin{equation*}
    \Delta b_k^{(m+1)}+8(N+1)^2r^{2N}e^{U_k}b_k^{(m+1)}=-E_k-\hat f(b_k^{(m)})
\end{equation*}
and $b_k^{(m+1)}(0)=0$.
Denote $F_k^m(r)=-E_k(r)-\hat{f}(b_k^{(m)}(r))$. Then by the iteration method as before if we set
$$\frac{d^2}{dr^2}b_k^{(m+1)}+\frac 1r\frac{d}{dr}b_k^{(m+1)}+8(N+1)^2r^{2N}e^{U_k}b_k^{(m+1)}=F_k^m, \quad 0<r<\tau \epsilon_k^{-1},$$
and 
$$b_k^{(m+1)}(0)=0.$$
The homogeneous equation has two fundamental solutions, one is 
$$u_1=\frac{1-r^{2N+2}}{1+r^{2N+2}}.$$
The second fundamental solution $u_2$ satisfies $|u_2(r)|\le C\log \frac{1}r$ near $0$ and $\infty$. We can construct $b_k^{(m+1)}$ from $b_k^{(m)}$ as
$$b_k^{(m+1)}(r)=u_1(r)\int_0^rtF_k^m(t)u_2(t){\rm d}t+u_2(r)\int_0^rtF_k^m(t)u_1(t){\rm d}t. $$
If $b_k^{(m)}(t)$ satisfies
$$|b_k^{(m)}(t)|\le C\epsilon_k^2\log (2+t),$$
one can verify by direct computation that $b_k^{(m+1)}$ satisfies the same bound. Thus by standard Brower fixed point theorem, there is a $b_k$ such that 
$$\Delta b_k+8(N+1)^2r^{2N}e^{U_k}b_k=-E_k-\hat{f}(b_k).$$
Theorem \ref{lap-non} is established. $\Box$

\end{document}